\documentclass[12pt,reqno]{article}

\usepackage[usenames]{color}
\usepackage{amssymb}
\usepackage{amsmath}
\usepackage{bbm}
\usepackage{amsthm}
\usepackage{amsfonts}
\usepackage{amscd}
\usepackage{bm}
\usepackage{graphicx}

\usepackage{stmaryrd}

\usepackage[colorlinks=true,
linkcolor=webgreen,
filecolor=webbrown,
citecolor=webgreen]{hyperref}

\definecolor{webgreen}{rgb}{0,.5,0}
\definecolor{webbrown}{rgb}{.6,0,0}

\usepackage{color}
\usepackage{fullpage}
\usepackage{float}

\usepackage{graphics}
\usepackage{latexsym}
\usepackage{epsf}
\usepackage{breakurl}

\usepackage{pgfplots}
\usepackage{filecontents}

\newcommand{\seqnum}[1]{\href{https://oeis.org/#1}{\underline{#1}}}

\usepackage[caption=false,font=footnotesize]{subfig}
\newcommand{\psubref}[1]{\protect\subref{#1}}

\begin{document}

\theoremstyle{plain}
\newtheorem{theorem}{Theorem}
\newtheorem{corollary}[theorem]{Corollary}
\newtheorem{lemma}[theorem]{Lemma}
\newtheorem{proposition}[theorem]{Proposition}

\theoremstyle{definition}
\newtheorem{definition}[theorem]{Definition}
\newtheorem{example}[theorem]{Example}
\newtheorem{conjecture}[theorem]{Conjecture}

\theoremstyle{remark}
\newtheorem{remark}[theorem]{Remark}

\begin{center}
\vskip 1cm{\LARGE\bf 
A new conjecture equivalent to\\ Collatz conjecture
}
\vskip 1cm
\large
Giulio Masetti\\
Istituto di Scienza e Tecnologia dell'Informazione,\\
Consiglio Nazionale delle Ricerche,\\
Via Moruzzi 1, Pisa, Italy,\\
Istituto Comprensivo di Anghiari e Monterchi,\\
Piazza del Popolo, Anghiari, Italy,\\
\href{mailto:giulio.masetti@isti.cnr.it}{\tt giulio.masetti@isti.cnr.it}
\end{center}

\vskip .2 in

\begin{abstract}
In this paper a new conjecture equivalent to Collatz conjecture is presented. 
In particural, showing that (all) the solution(s)
of newly introduced iterative functional equation(s) have a given property
is equivalent to prove Collatz conjecture.
\end{abstract}

\section{Introduction}\label{sec:intro}
In the following a new conjecture equivalent to Collatz conjecture will be presented.
Calling $\mathbb N=\{1,2,\dots\}$ the set of natural numbers, consider the map 
$\mathrm{Col}:\mathbb N\rightarrow\mathbb N$ defined by
\begin{equation}
\label{col}
\mathrm{Col}(n) := 
\begin{cases}
	3n+1, &\text{ if } n\text{ is odd};\\
	\frac{n}{2}, &\text{ if } n\text{ is even},	
\end{cases}
\end{equation}
and, given $k\in\mathbb N$, the sequence $(d_{k,i})_{i\ge 0}$ defined by
\begin{equation}
\label{eq:departure}
d_{k,i} :=
\begin{cases}
	k, &\text{ if } i=0;\\
	\mathrm{Col}(d_{k,i-1}), &\text{ if } i>0.
\end{cases}
\end{equation}
One of the formulations of the Collatz conjecture can be stated as:
\begin{conjecture}[Collatz conjecture]
	\label{conj:Collatz}
	For all $k\ge 1$, there exists $i\ge0$ such that $d_{k,i}=1$.
\end{conjecture}
Several attempts have been made for proving or disproving 
conjecture~\ref{conj:Collatz}, none reached the goal.
For $k=1,2,3$ and $4$ the \emph{departure sequence} defined in eq.~\ref{eq:departure},
also called the $3n+1$ sequence, has a simple behavior, but for $k\ge 5$ the behavior is
erratic.  

\section{From Departure to Arrival}
\label{sec:idea}
Eq.~\ref{eq:departure} defines a discrete dynamical system,
i.e., $i$ is usually interpreted as a discrete time instant,
$d_{k,i}$ spans the states of the system and the system behavior 
is characterized answering several (infinite) times to the question: 
``once the system is in this state, where it can go?''.
The focus is on \emph{departure}.

In this paper, instead, the focus is on \emph{arrival}, 
examining the question ``from which states
the system can arrive in a given state?''.
The shift of perspective is supposed to hide unnecessary details
because conjecture \ref{conj:Collatz} is formulated in terms
of arrival (will $(d_{k,i})_{i\ge 0}$ ever arrive to $1$?)
and then packing the available information prioritizing the receiver
instead of the sender restates the problem in the language of the expected solution.
The agenda is:
\begin{enumerate}
	\item following a classical approach, introduce invariant measures
	with respect to the dynamical system. Every specific instance of 
	conjecture \ref{conj:Collatz} corresponds to find an invariant measure,
	and this can be stated as finding one solution of a quite general and elegant
	functional equation (the solution being a probability mass function associated with the measure);
	\item define the new sequence $(a_{k,n})_{n>0}$, from now on called
	\emph{arrival sequence}, in eq.~ \ref{eq:arrival}.
	The idea is that the arrival sequence marks if $n$ is a state of
	the system, i.e., if $n$ is reachable starting from $k$ through a certain number
	of applications of $\mathrm{Col}$.
	More formally:
	\begin{equation}
	\label{eq:implies}
	a_{k,n}\neq0 \implies \exists i\ge 0\text{ such that }d_i=n;
	\end{equation}
	\item following the generatingfunctionology approach \cite{gf}, 
	define the formal power series $A_k(x):=\sum_{n\ge 0}a_{k,n} x^n$
	and find an equation that it solves (i.e., eq.~\ref{fps}).
	This way it is possible to formulate a new conjecture (conjecture \ref{conj:arrival})
	that is equivalent to conjecture \ref{conj:Collatz}. Notice that $a_{k,n}$
	stores interesting information linked to the state hitting time, 
	i.e., how many times the dynamical system touches $n$, 
	but for the scope of this paper it is relevant just to know if $a_{k,n}$ is nonzero;
	\item evaluate $A_k(x)$ in $x^2$ and, by equating the (formal) derivatives of left 
	and right hand sides of eq.~\ref{fps_square}, verify the expected behavior 
	and prove useful properties of the arrival sequence;
	\item proving that, for all $k\ge 1$,
	$A_k'(0)=a_1$ is nonzero means proving conjecture \ref{conj:arrival}.
	This step is out of the paper scope (i.e., at the moment, there is no proof of conjecture \ref{conj:arrival}), few attempts are sketched in Section \ref{sec:attempts}.
\end{enumerate}

Switching from departure to arrival can also shead new light on the
relationship among the $3n+1$ problem, the $3n-1$ problem and
variants of the aforementioned, such as
\begin{equation}
	\label{eq:generalization}
	T_h(n) := 
	\begin{cases}
		h\cdot n\pm 1, &\text{ if } n\text{ is odd};\\
		\frac{n}{2}, &\text{ if } n\text{ is even},	
	\end{cases}
\end{equation}
with $h$ odd.

\section{$\mathrm{Col}$-invariant measures over $\mathbb N$}
\label{sec:arrival}
A measure $\mu$ over $\mathbb N$ is $\mathrm{Col}$-invariant if,
for every measurable set $E$, we have $\mu(E)=\mu\big(\mathrm{Col}^{-1}(E)\big)$.
In the following we will focus on measures $\mu$ such that:
\begin{itemize}
	\item for all $n\in\mathbb N$, the set $\{n\}$ is $\mu$-measurable,
	and the measure mass function $a_n:=\mu\big(\{n\}\big)$ defines
	the \emph{arrival} sequence $(a_n)_{n>0}$,
	\item $\mu$ is (finitely) additive.
\end{itemize}

The departure sequence, depicted in Figure \ref{fig:3n+1}, 
clearly identifies local relations among few states of the system, 
so the arrival sequence defines the recurrence relation: 
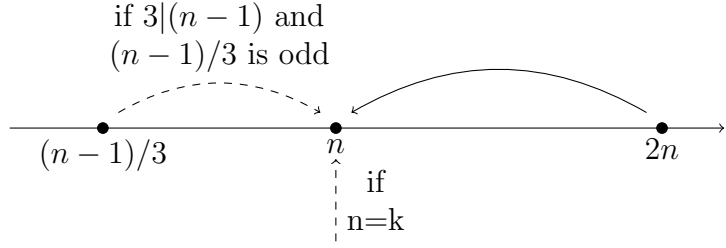
\begin{figure}
	\centering
	\begin{tikzpicture}
	% n=2, 7 and 14 are scaled by 1.5
	\draw[->] (0,0) -- (23em,0);
	\filldraw[black] (3em,0) circle (2pt) node[below] (f) {$(n-1)/3$};
	\filldraw[black] (10.5em,0) circle (2pt) node[below] (n) {$n$};
	\filldraw[black] (21em,0) circle (2pt) node[below] (b) {$2n$};
	\path[->,dashed] (3.5em,.5em) edge [bend left] node[above,align=center] 
	{if $3|(n-1)$ and\\$(n-1)/3$ is odd} (10em,.5em);
	\path[->] (20.5em,.5em) edge [bend right] node[above] {} (11em,.5em);
	\node at (10.5em, -4em) (s) {};
	\draw[->,dashed] (s) -- node[right,align=center] {if\\n=k} (10.5em,-1em);
	\end{tikzpicture}
	\caption{Pictorial representation of $\mathrm{Col}$ from the point
		of view of the receiver instead of the sender. Dashed arrows represent conditional steps.}
	\label{fig:3n+1}
\end{figure} 
\begin{equation}
\label{eq:arrival_measure}
a_{n} := a_{2n}
+\mathbbm{1}_{3|(n-1)\text{ and } \frac{n-1}{3}\text{ is odd}}
\cdot a_{\frac{n-1}{3}},
\end{equation}
Notice that there are infinite measures $\mu$ for which $a_k>0$.
It is now possible to formulate a new conjecture, equivalent to
conjecture \ref{conj:Collatz}:
\begin{conjecture}[arrival conjecture]
	\label{conj:arrival_measure}
	For all $k\ge 1$ such that $a_k>0$ we have $a_{1}> 0$.
\end{conjecture}

In the following we will focus mainly on another recurrence relation,
the one that corresponds to the normalization $a_k:=1$, and
we will call also the corresponding sequence the arrival sequence:
\begin{equation}
	\label{eq:arrival}
	a_{k,n} := \delta_{n-k} + \cdot a_{k,2n}
	+\mathbbm{1}_{3|(n-1)\text{ and } \frac{n-1}{3}\text{ is odd}}
	\cdot a_{k,\frac{n-1}{3}},
\end{equation}
It is possible then to define a new conjecture that is equivalent to
conjecture \ref{conj:Collatz}:
\begin{conjecture}[arrival conjecture]
	\label{conj:arrival}
	For all $k\ge 1$, $a_{k,1}\neq 0$.
\end{conjecture}

Notice that whenever the departure sequence touches 
$n=2^{h}\cdot o$, where $o$ is an odd number, the arrival sequence
marks as nonzero $a_{2^i\cdot o}$ for $0\le i\le h$.
In other words, the nonzero elements of the arrival sequence are grouped in subsets, 
each subset being labeled by an odd number.

Defining the probability mass function 
$A(x):=\sum_{n\ge 0}a_{n} x^n$
and the formal power series $A_k(x):=\sum_{n\ge 0}a_{k,n} x^n$, 
the first addend of eq.~\ref{eq:arrival_measure}
produces $x^k$, the second addend in eq.~\ref{eq:arrival}, 
from now on called \emph{backward} term,
produces $(A_k(\sqrt{x})+A_k(-\sqrt{x}))/2$ 
and the third addend, called \emph{forward} term, produces
$x\big(A_k(x^3)-A_k(-x^3)\big)/2$.
By grouping all the terms we obtain
\begin{align}
A(x) &= \frac{A\big(\sqrt{x}\big)+A\big(-\sqrt{x}\big)}{2}
+ \frac{x\big(A\big(x^3\big)-A\big(-x^3\big)\big)}{2},\label{fps_measure}\\
A_k(x) &= \frac{A_k\big(\sqrt{x}\big)+A_k\big(-\sqrt{x}\big)}{2}
+ \frac{x\big(A_k\big(x^3\big)-A_k\big(-x^3\big)\big)}{2} + x^k.\label{fps}
\end{align}

Actually, we can slightly modify $A_k$ introducing weights
in front of the backward and forward terms carefully chosen to make
it analytic (and then eq.~\ref{fps} a functional equation):
\begin{equation}
	\label{fun_eq}
	A_k(x) := w_b\frac{A\big(\sqrt{x}\big)+A\big(-\sqrt{x}\big)}{2}
	+ w_f\frac{x\big(A\big(x^3\big)-A\big(-x^3\big)\big)}{2} + x^k.
\end{equation}
In the rest of the paper $A_k$ will indicate this modified expression.
We can also evaluate eq.~\ref{fun_eq} in $x^2$ for being able to differentiate 
\begin{equation}
\label{fps_square}
2A_k(x^2) = w_bA_k(x) + w_b A_k(-x) +w_f x^2A_k(x^6)-w_fx^2A_k(-x^6) +2x^{2k}.
\end{equation}
Notice that working with (formal) power series makes it possible to
ignore the dynamics' details and offers a global view on the problem at hand.
In particular, conjecture \ref{conj:Collatz} is equivalent to:
\begin{conjecture}[formal power series conjecture]
	\label{conj:arrival_fps}
	For all $k\ge 1$, $A_k'(0)=a_{k,1}\neq 0$.
\end{conjecture}

Summing up: $A(x)$ in eq.~\ref{fps_measure} captures
all the relevant information of $\mathrm{Col}$;
$A_k(x)$ in eq.~\ref{fun_eq} caputres a specific dynamics, the one
that touches $k$.

Following the same reasoning, eq.~\ref{eq:generalization}
defines
\begin{equation*}
	A_k(x) := w_b\frac{A\big(\sqrt{x}\big)+A\big(-\sqrt{x}\big)}{2}
	+ w_f\frac{x^{\pm 1}\big(A\big(x^h\big)-A\big(-x^h\big)\big)}{2} + x^k.
\end{equation*}

\section{Related Work}
\label{related}
An account of several strategies to address the Collatz conjecture can be 
found in papers of Lagaris \cite{Lagaris03,Lagaris06} and citations therein.

Notice that $A(x)$ is a link between the Collatz conjecture and ergodic theory, but with
a different flawor with respect to the link discussed by Lagaris (Section 2.8 of \cite{Lagaris85}).
In addition, the (functional) equations presented in this paper 
are quite different from the those presented by Berg and Burckel \cite{Berg,Burckel94},
but quite similar to the one presented by Neklyudov \cite{Neklyudov21}.  

$A(x)$ and $A_k(x)$ present some similarity with 
the functional equations addressed in works of Heuvers, Moak, Boursaw, Odlyzko and Hardy \cite{SquareRootSpiral,Odlyzko82,Hardy}.

\section{Investigate numerically the arrival sequence}
\label{sec:verification}
Equating\footnote{To speed up computations, at \url{https://github.com/106ohm/CollatzConjectureFunctionalEquation} 
a simple Maxima program is available.} 
derivatives of even order (those of odd order always produces $0=0$)
of left and right hand sides of eq.~\ref{fps_square}, an infinite set of linear equations
can be defined. For instance, fixing $k=5$ the first $32$ derivatives produce:
\begin{align*}
a_{5,1} &= w_b a_{5,2} & a_{5,2} &= w_b a_{5,4} & a_{5,3} &= w_b a_{5,6} & a_{5,4} &= w_b a_{5,8} + w_f a_{5,1}\\
a_{5,5} &= w_b a_{5,10} + 1 & a_{5,6} &= w_b a_{5,12} & a_{5,7} &= w_f a_{5,2} + w_b a_{5,14} & a_8 &= w_b a_{5,16}\\
a_{5,9} &= w_b a_{5,18} & a_{5,10} &= w_f a_{5,3} + w_b a_{5,20} & a_{5,11} &= w_b a_{5,22} & a_{5,12} &= w_b a_{5,24}\\
a_{5,13} &= w_f a_{5,4} + w_b a_{5,26} & a_{5,14} &= w_b a_{5,28} & a_{5,15} &= w_b a_{5,30} & a_{5,16} &= w_f a_{5,5} + w_b a_{5,32}
\end{align*}
Thus, in general, $(1-w_b^2w_f)a_1=w_b^3a_{5,8}$, and for $k=5$, $a_{5,8}=w_bw_f a_{5,5}+w_b^2 a_{5,32}$.
Being $a_{5,5}\neq 0$, if
\begin{equation}
\label{eq:condition_w}
w_b^2 w_f\neq 1
\end{equation}
then $a_{5,1}\neq 0$, so conjecture \ref{conj:arrival_fps} is verified for $k=5$.  

Rearranging eq.~\ref{fun_eq}, 
and following the general approach of de Bruijn \cite{DeBruijn79}, 
it is easy to define $\{A_{k,i}\}_{i=0,1,\dots}$ such that
$A_k(x)=\lim_{i\rightarrow\infty}A_{k,i}(x)$.
Calling $g_k(x)=x^k$, $\varphi(x)=\sqrt{x}$ and $\psi(x)=x^3$, we have
\begin{align}
	&F(y,x) = \frac{w_b}{2} \Big(y(\varphi(x))+y\big(-\varphi(x)\big)\Big) 
	+ \frac{w_f}{2} x \Big(y\big(\psi(x)\big)-y\big(-\psi(x)\big)\Big) + g_k(x),\label{eq:functional}\\
	&\begin{cases}
		A_{k,0}(x) &= g_k(x),\\
		A_{k,i+1} &= F(A_{k,i}(x),x),\text{ for } i=0,1,\dots, 
	\end{cases} \label{eq:iterate}
\end{align}
and then, fixed reasonable values of $w_b$ and $w_f$, we can plot $A_{k,i}$ for a given $k$.
\begin{figure}
	\centering
	\subfloat[$A_{5,4}(x)$]{
		\includegraphics[scale=.25]{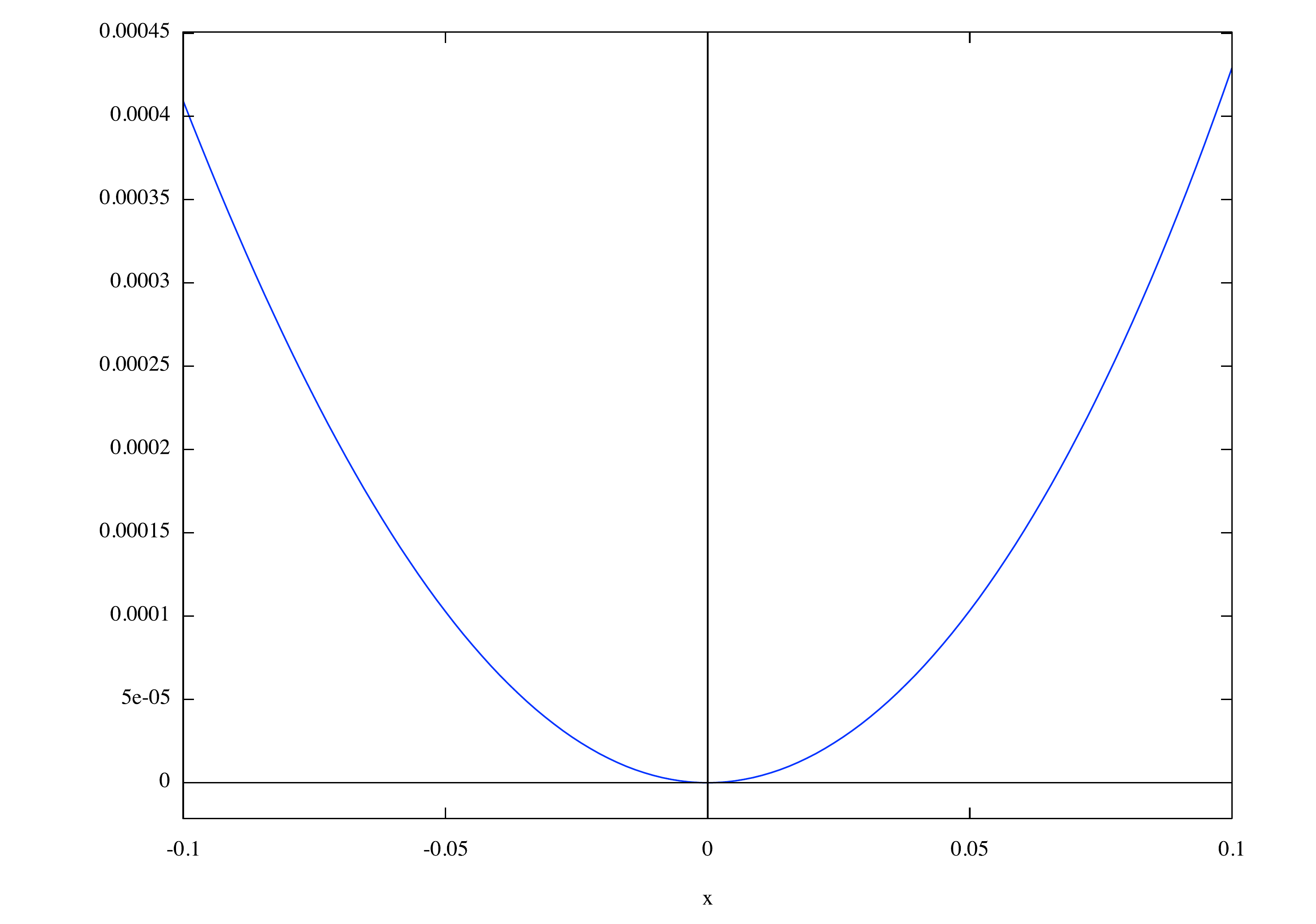}
		\label{fig:k5iterate4small}
	}
	\subfloat[$A_{5,5}(x)$]{
		\includegraphics[scale=.25]{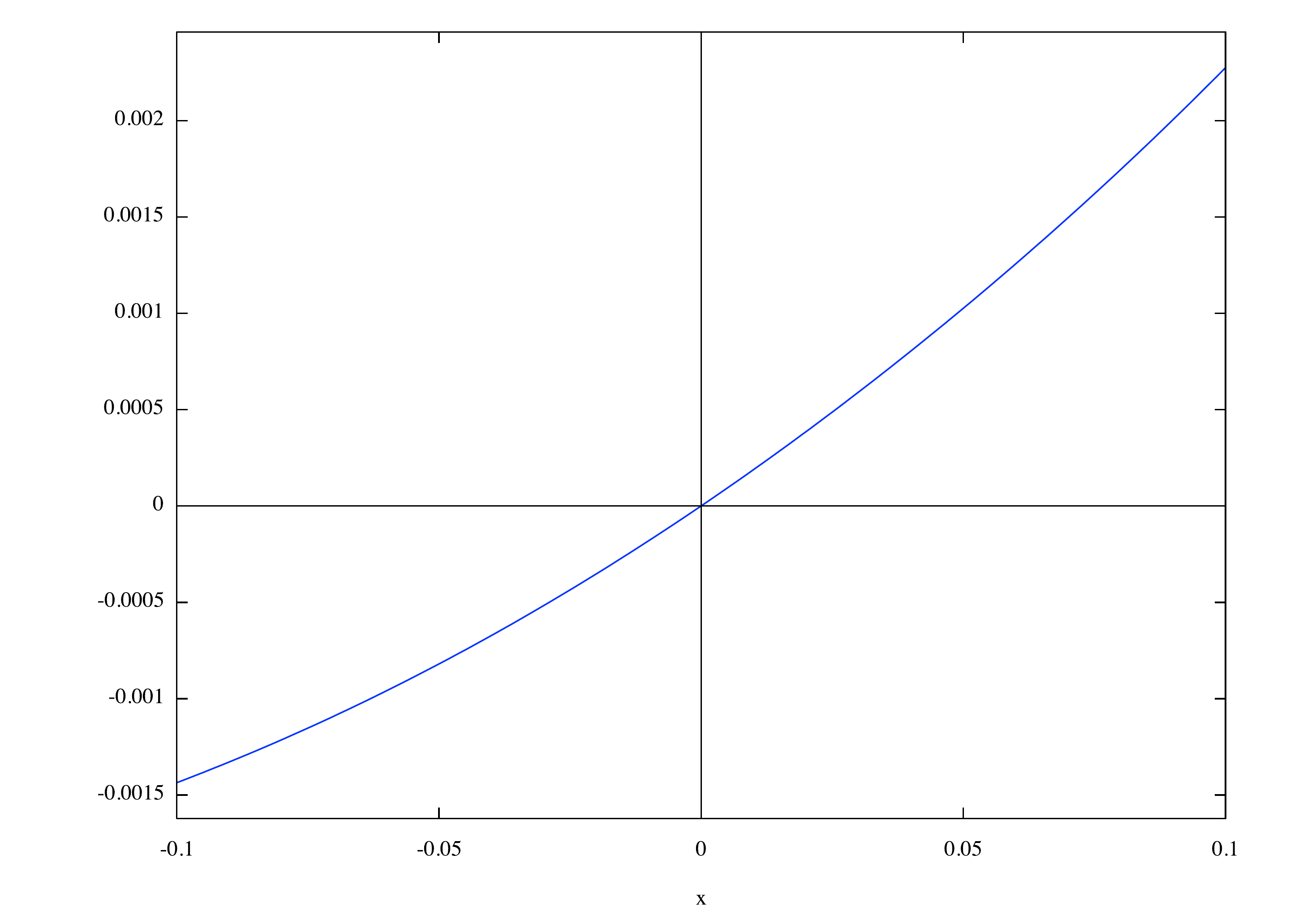}
		\label{fig:k5iterate5small}
	}
	\\
	\subfloat[$A_{5,5}(x)$]{
		\includegraphics[scale=.25]{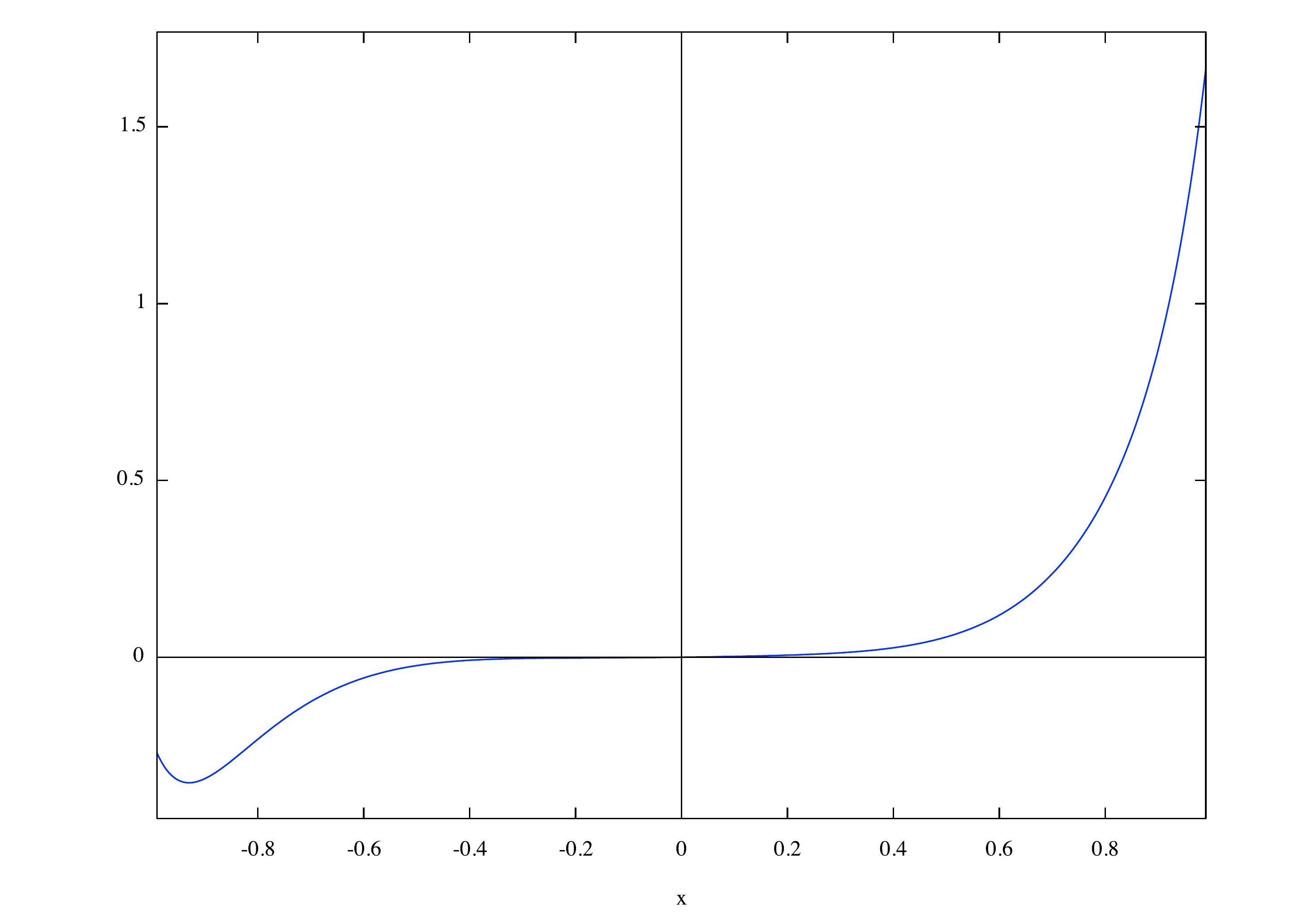}
		\label{fig:k5iterate5}
	}
	\subfloat[$A'_{5,5}(x)$]{
		\includegraphics[scale=.25]{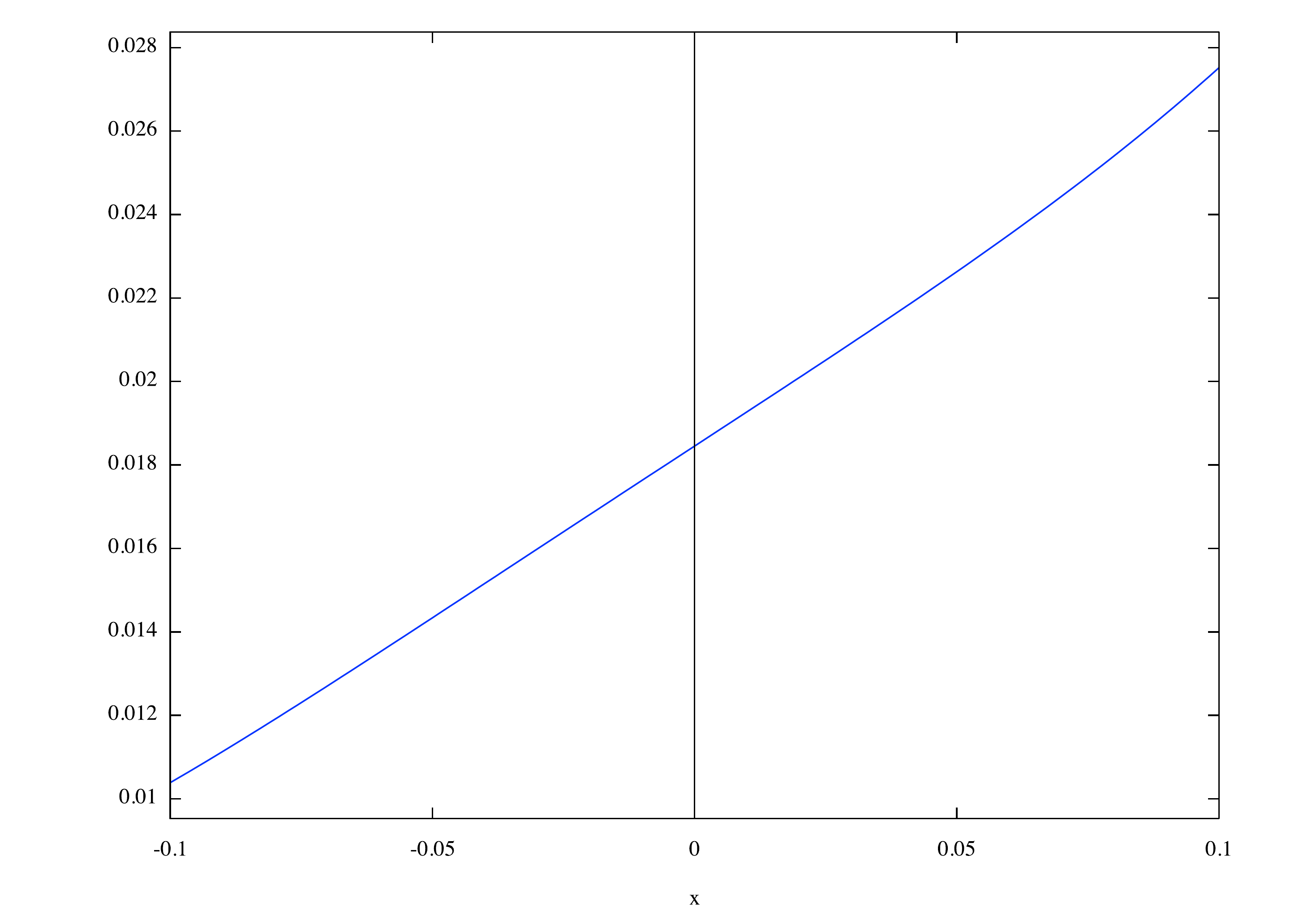}
		\label{fig:k5iterate5derivative}
	}
	\label{fig:k5_iterate}
	\caption{Study of $A_{5,i}$ as a real function, i.e., $x\in\mathbb R$ and $A_{5,i}(x)\in\mathbb R$. Parameters setting: $k=5$, $w_b=0.45$ and $w_f=0.45$ for $A_{5,4}$ 
		in \psubref{fig:k5iterate4small} 
		and $A_{5,5}$ in \psubref{fig:k5iterate5small} on $x\in[-0.1,0.1]$.
		The function $A_{5,5}$ on $x\in[-0.99,0.99]$ is plotted in \psubref{fig:k5iterate5},
		and $A'_{5,7}$ in \psubref{fig:k5iterate5derivative} to check that $A'_{5,5}(0)\neq 0$.}	
\end{figure} 
For instance\footnote{Code available at \url{https://github.com/106ohm/CollatzConjectureFunctionalEquation}.},
in Figure \ref{fig:k5_iterate}, where $k=5$, it is clear that for
$i\le 5$ the value of $A'_{k,i}(0)$ is zero, but for $i=6$ there
is a qualitative change: $A'_{k,i}(0)\neq 0$.

Actually, $A_{k,i}$ mimics the departure sequence in eq.~\ref{eq:departure} in the sense
that relevant information is accumulated in $A_{k,i}$ at increasing of $i$.
For instance, Figure \ref{fig:k5_iterate} shows the qualitative change in $A_{k,i}$
switching from $i\le 5$ to $i\ge 6$, the same change that we see in the departure sequence.

\begin{figure}
	\centering
	\begin{tikzpicture}
		\draw[->] (-4em,0em) -- (4em,0em);
		\draw (-3em,-3em) -- (3em,3em);
		\draw[->] (0em,-4em) -- (0em,4em);
		\draw (-3em,3em) -- (3em,-3em);
		\node[fill, draw, circle, minimum size=2, inner sep=1pt, label= below:$0.2$] (1) at (3em,0em) {};
		\node[fill, draw, circle, minimum size=2, inner sep=1pt] (2) at (3em,3em) {};
		\node[fill, draw, circle, minimum size=2, inner sep=1pt, label= left:$0.2i$] (3) at (0em,3em) {};
		\node[fill, draw, circle, minimum size=2, inner sep=1pt] (4) at (-3em,3em) {};
		\node[fill, draw, circle, minimum size=2, inner sep=1pt] (5) at (-3em,0em) {};
		\node[fill, draw, circle, minimum size=2, inner sep=1pt] (6) at (-3em,-3em) {};
		\node[fill, draw, circle, minimum size=2, inner sep=1pt] (7) at (0em,-3em) {};
		\node[fill, draw, circle, minimum size=2, inner sep=1pt] (8) at (3em,-3em) {};
		
		\draw[>>>->>>] (5) -- (1);
		\draw[>>>->>>] (7) -- (3);
		\draw[>>>->>>] (4) -- (8);
		\draw[>>>->>>] (6) -- (2);
		%%%
		\node (0) at (14em,0em) {\includegraphics[scale=0.2]{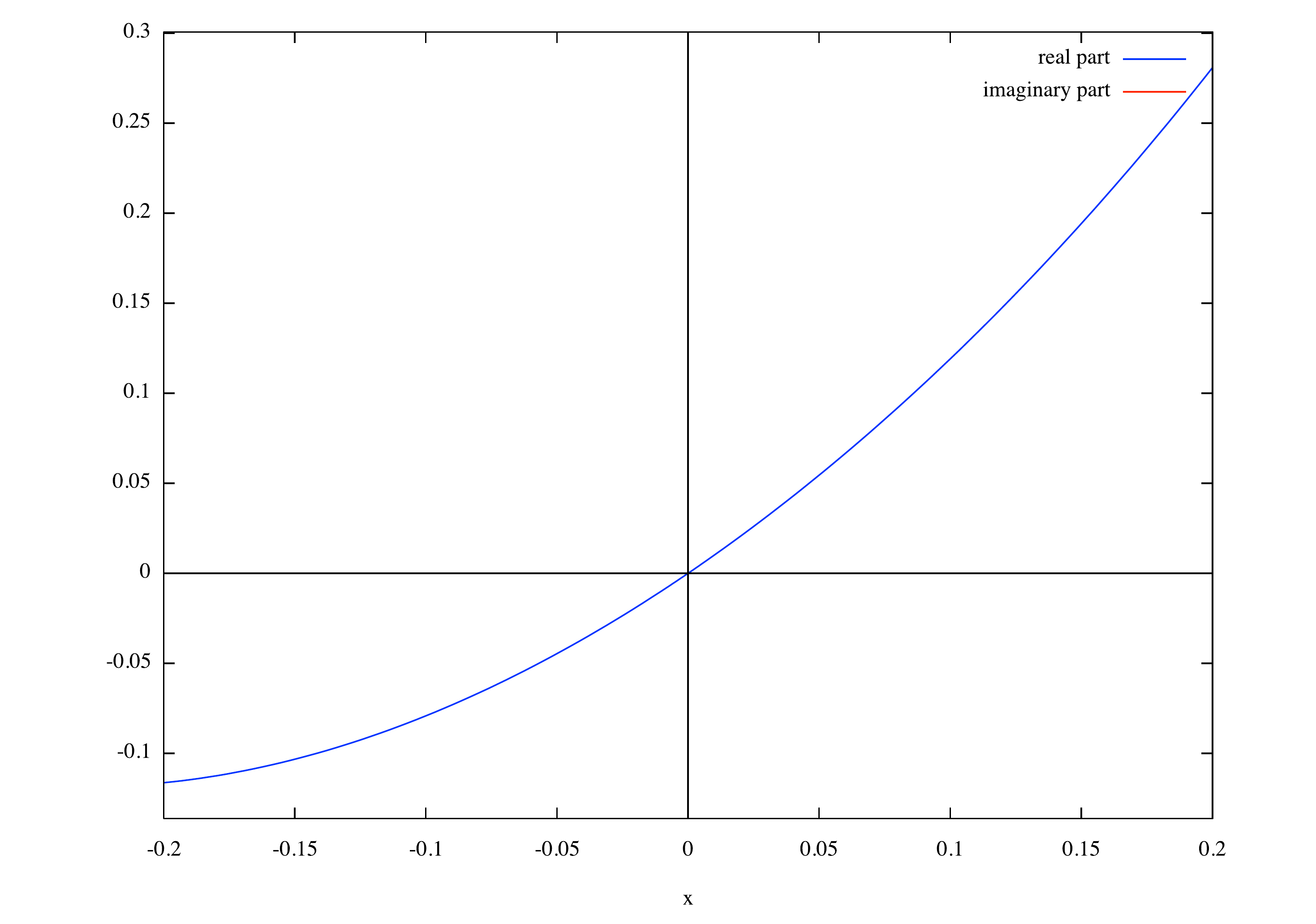}};
		\node (1) at (14em,12em) {\includegraphics[scale=0.2]{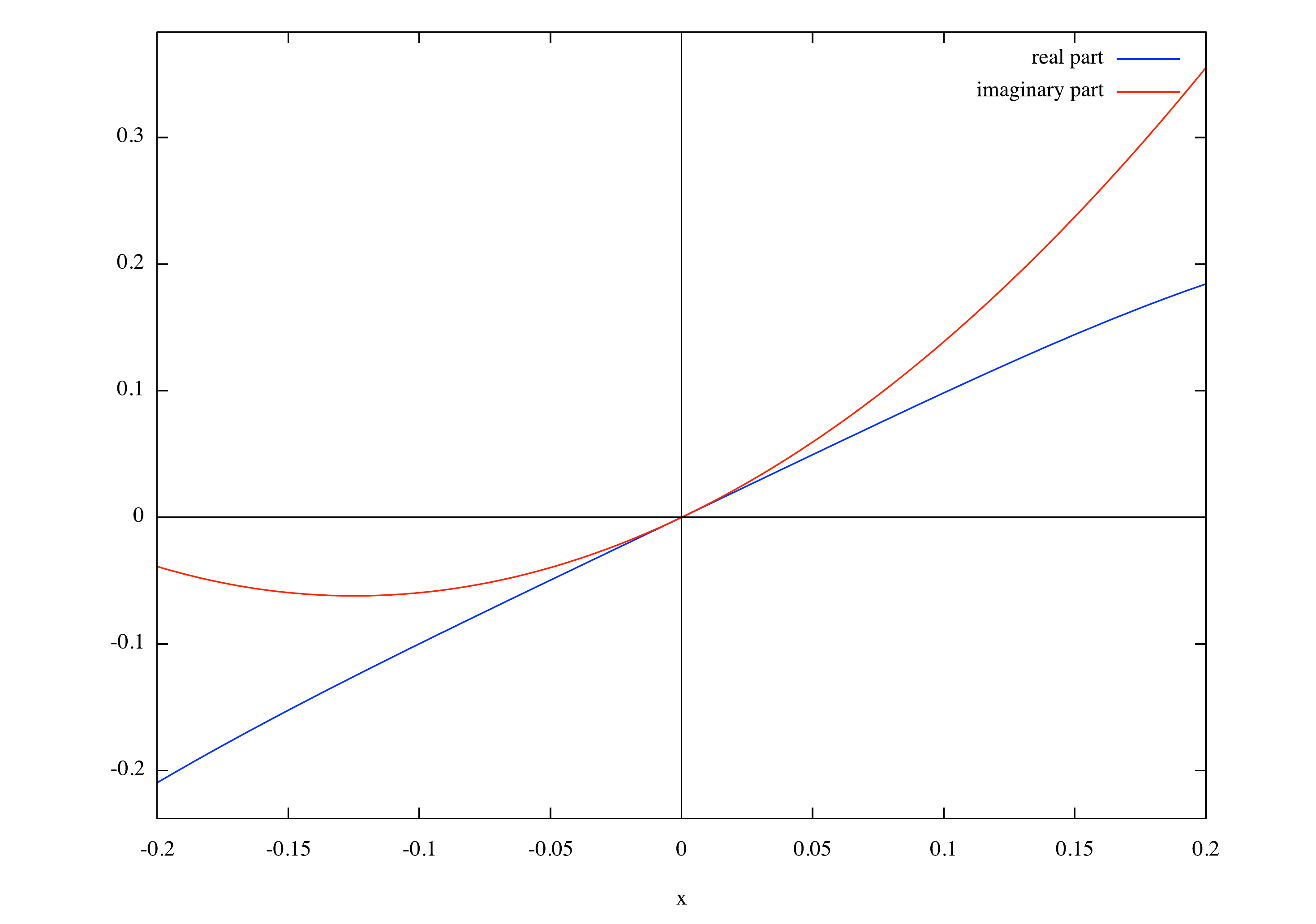}};
		\node (2) at (0em,12em) {\includegraphics[scale=0.2]{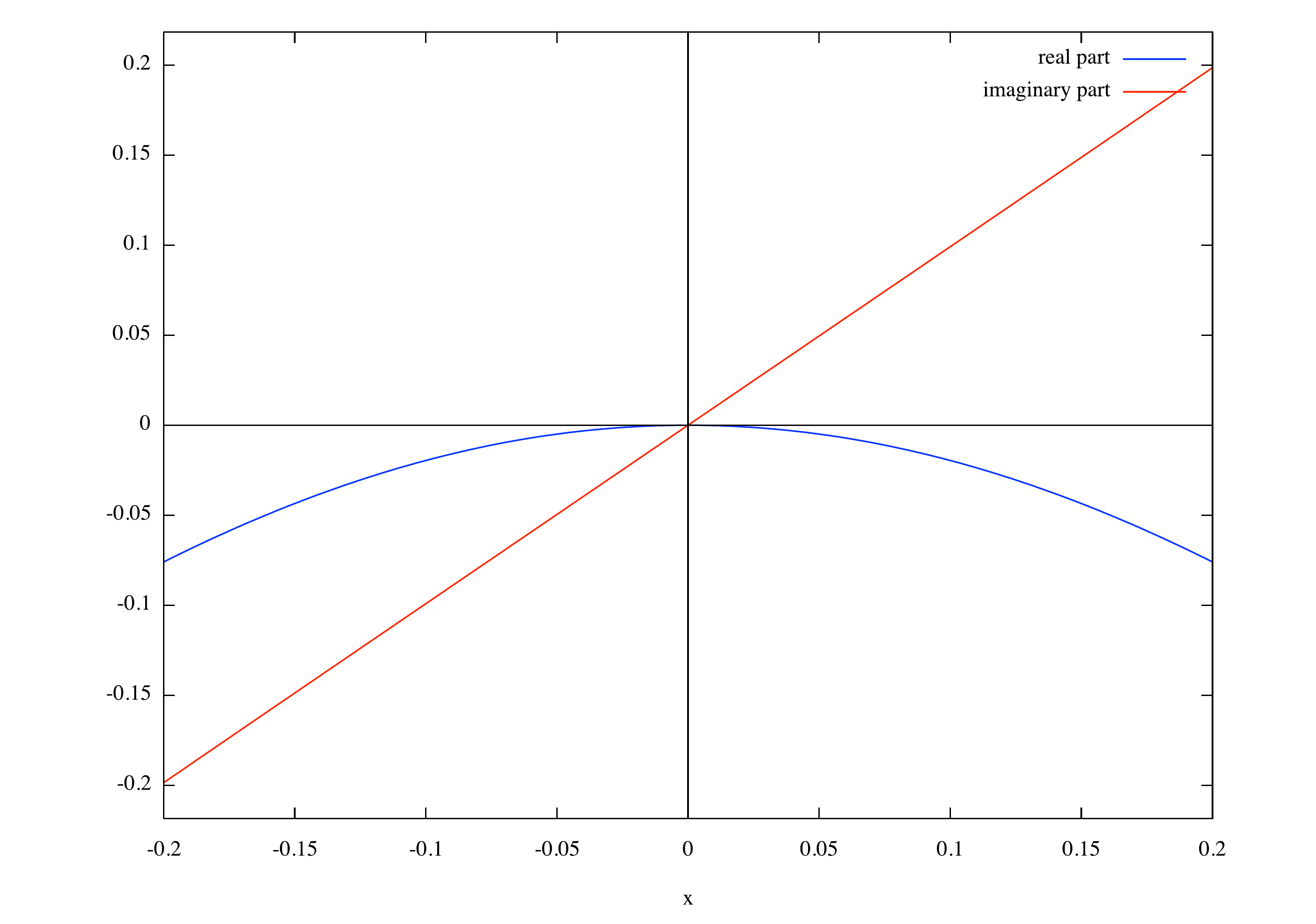}};
		\node (3) at (-14em,12em) {\includegraphics[scale=0.2]{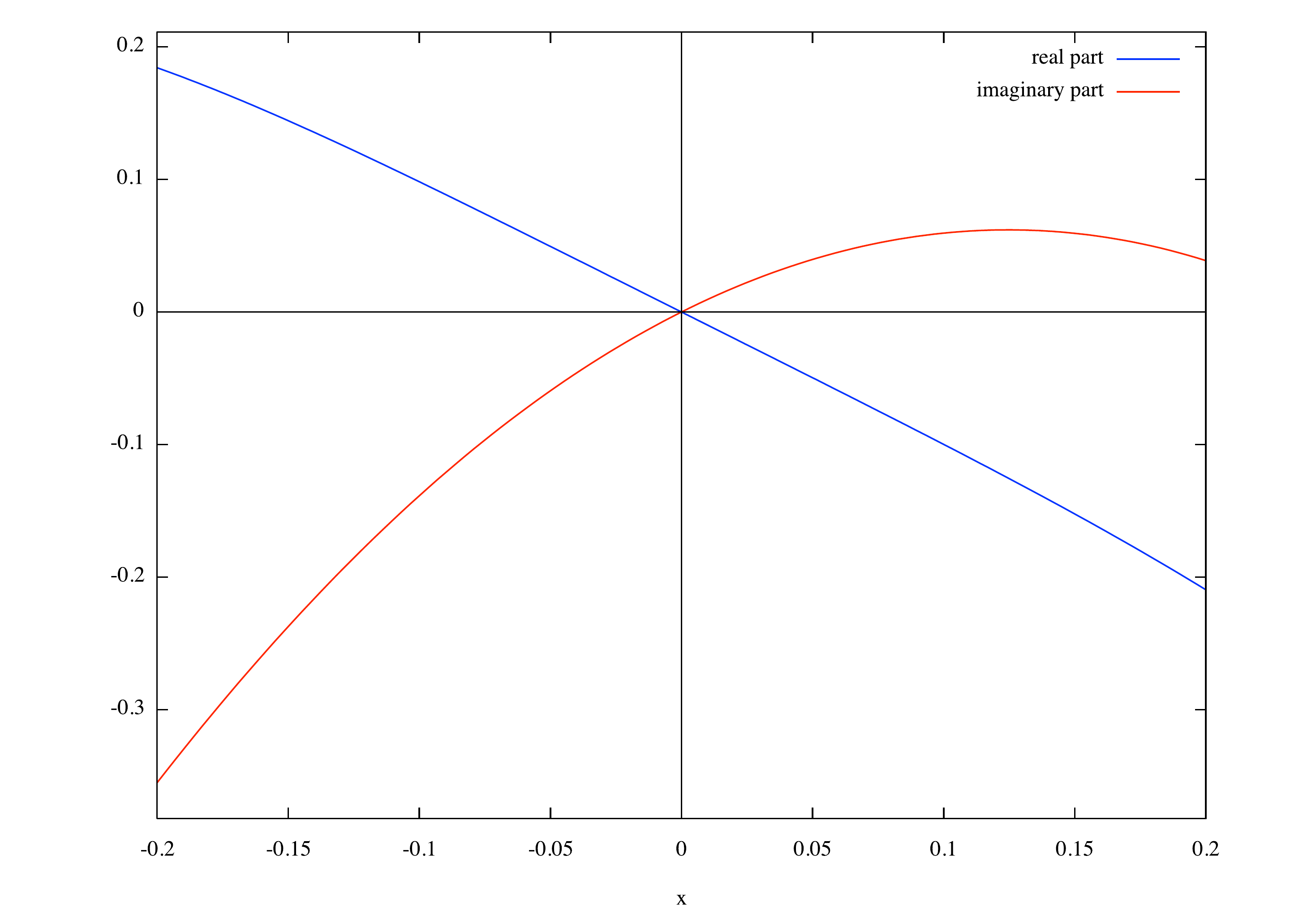}};
		\draw (2em,0em) edge[bend left,->] node[above] {$A_{5,7}(x)$} (0.west);
		\draw (2em,2em) edge[bend left,->] node[below right] {$A_{5,7}(x)$} (1.south west);
		\draw (0em,2em) edge[bend left,->] node[above left] {$A_{5,7}(x)$} (2.south);
		\draw (-2em,2em) edge[bend left,->] node[below] {$A_{5,7}(x)$} (3.south);
	\end{tikzpicture}
	\label{A7}
	\caption{Visualizing how $A_{5,7}(x)$ behave when $x$ crosses 0, i.e., 
	for $x$ on the marked oriented segments with the following parameters: 
	$k=5$, $w_b=0.999$ and $w_f=0.995$.
	Of course, on the real line $A_{5,7}(x)\in\mathbb R$, so the imaginary part is zero.
	}
\end{figure}

Notice that, clearly, $A_{k,i}$ is the sum of $i+1$ terms involving the nested evaluation of
$g$, $\varphi$ and $\psi$. 
Unfortunately, it is not as ``easy'' as for the functional equations introduced by Odlyzko or Hardy \cite{Odlyzko82,Hardy} to directly check that $A^{\prime}_{k,i}(0)\neq 0$ for $i$ large enough.

\section{$A_k(x)$ is analytic under reasonable assumptions}
\label{sec:analytic}
Eq.~\ref{fun_eq} is, following the classification presented 
by Acz\'el \cite{Aczel1984FunctionalEH}, a linear iterative functional equations.
As such, it is quite easy to prove that its solution is unique and continuous
(see Theorem and Lemma 6.6.1 developed by Kuczma \cite{Kuczma90}) when $|w_b|,|w_f|<1$.
Further, following a similar argument, it is possoble to show that
$A_k(x)$ is analytic.
Indeed, 
focusing on the space of bounded analytic functions on $B(0,r)$, 
we can define the operators:
\begin{align}
\big(H_1g\big)(x) &:= \frac{1}{2}\sum_{z^2=x} g(z),\label{H1}\\
\big(H_2g\big)(x) &:= \frac{x}{2}\big(g\big(x^3\big)-g\big(-x^3\big)\big),\label{H2}\\
\big(Hg\big)(x) &:= w_b H_1 + w_f H_2.\label{H} 
\end{align}
and observe that
$H_1$ and $H_2$ are linear and both bounded for $r=1$: $||H_1||<1$ and a rough estimate
for $H_2$ is $||H_2||<r=1$.
Thus, if $|w_b|+|w_f|<1$ then $||H||<1$ and $(I-H)^{-1}=\sum_{i\ge 0}H^i$ makes sense;
so we can rewrite eq.~\ref{fun_eq} as
\begin{equation}
\label{eq:linear_operator}
A_k(x) = \big(I-H\big)^{-1}x^k = \sum_{n\ge 0} H^n x^k.
\end{equation}
One issue with eq.~\ref{eq:linear_operator} is that $H_1$ and $H_2$ does not commute
and then we cannot exploit the binomial theorem 
for $(w_b H_1+w_f H_2)^n$.
Nevertheless, $A(x)$ is analytic.

%Notice that $H^2_2=0$ and $H_1H_2 \big(f(x)\big) = H_2\big(f(\sqrt{x})\big)$, thus
%$H_2^{m_2}H_1^{m_1}H_2=0$ for $m_2\ge 1$ and $m_1\ge 0$.
%This means (\textcolor{red}{is this correct?}) that
%\[
%(w_b H_1+w_f H_2)^m \big(f(x)\big) = w_b^mH_1^m\big(f(x)\big)
%+w_fw_b^mH_2\big(f(x^{\frac{1}{2^m}})\big)
%+w_fw_b^m(H_2H_1^m)\big(f(x)\big).
%\]
%In particular, 
%\begin{align}
%A(x) &= \sum_{m\ge 0} \Big(w_b^mH_1^m\big(x^k\big)
%+w_fw_b^mH_2\big(x^{\frac{k}{2^m}}\big)
%+w_fw_b^m(H_2H_1^m)\big(x^k\big)\Big)\nonumber\\
%&=\sum_{m\ge 0,z^{2^m}=x} (w_b/2)^m z^k 
%+\sum_{m\ge 0}w_fw_b^mH_2\big(x^{\frac{k}{2^m}}\big)
%+\sum_{m\ge 0,z^{2^m}=x}w_f(w_b/2)^m H_2\big(z^k\big).\label{solution}
%\end{align}
%
%\textcolor{red}{Is the term $H_2H_1$ problematic?} I mean:
%the square root of $-x^3$ require to exploit the principal branch Log\dots 

Notice that, being $A_k(x)$ absolutely convergent and $a_n\in\mathbb R$
because $w_b,w_f\in\mathbb R$, if $x\in\mathbb R$ then $A(x)\in\mathbb R$.

\section{Trying to prove conjecture \ref{conj:arrival_fps}}
\label{sec:attempts}
The following list of trials is instructive by itself and, more importantly,
shows how conjecture \ref{conj:arrival_fps} is linked to other well known
topics or issues in (complex) analysis:
\begin{itemize}
	\item try to solve eq.~\ref{fun_eq} finding
	a closed formula (integral form, infinite products, etc),
	\item exploit Lagrange transform. Unfortunately the substitution
	of interest is the one with $b=1/2$ that is an open problem \cite{Buschman},
	\item exploit Mellin transform. Unfortunately, to the best of the author knowledge, 
	$\mathcal M\{f(\alpha\cdot t);s\}$ makes sense only for $\alpha>0$ \cite{TablesIntegralTransforms},
	\item define the Dirichlet series $D_k(s):=\sum_{n\ge 1}a_n/n^s$ 
	instead of the power series in eq.~\ref{fun_eq} and
	write the corresponding functional equation\footnote{Through the comparison with the Riemann zeta function
	it is easy to see that $D_k(x)$ is actualy an holomorphic function,
	being $0\le a_n<1$},
	\item try to exploit standard techniques for solving functional equations \cite{Aczel1984FunctionalEH,Aczel64}, 
	\item try to adapt the reasoning developed by Heuvers, Moak and Boursaw \cite{SquareRootSpiral} for instance
	exploiting the principal logarithm. %$y:=\mathrm{Log}(x)$
%	with $-\pi < \theta\le \pi$, it is possible to obtain the functional equation
%	\[
%	B_k(y) = \frac{w_b}{2}\big(B_k(y/2)+B_k(i\pi + y/2)\big)
%	+ \frac{w_f e^y}{2}\big(B_k(3y)-B_k(i\pi +3y)\big) + e^{ky}
%	\] 
	Unfortunatelly it is not easy to exploit the solution of the Kr\"ul equation \cite{SquareRootSpiral},
	\item the kind of manipulations developed by Hardy, de Bruijn and Odlyzko \cite{Hardy,DeBruijn79,Odlyzko82} 
	are also appealing, but (at the moment) not productive,
	\item finding a closed formula for $A_{k,i}(x)$ of eq.~\ref{eq:functional}, not assuming that is a polynomials, is complicated, 
	\item try to prove that $A$ is injective in a neighborhood of $0$, thus $A'(0)\neq 0$. 
	This follows from a nontrivial result of complex analysis
	(see Theorem 7.4 of Conway's book \cite{Conway}).
\end{itemize}

\section{Acknowledgments}\label{sec:ack} 
Many are the people I asked suggestions on parts of the reasoning,
often without mentioning them my aim. If you find mistakes in my developments
please consider them as my faults, not their. 
In particular, I would like to thank Pietro Majer, Edward Kamen, Leonardo Robol, 
Francesca Acquistapace, Fabrizio Broglia and Felicita Di Giandomenico for their help, 
and my wife Francesca Papini for her support during a complex period of our life.

\bibliographystyle{jis} 
\bibliography{main}

\bigskip
\hrule
\bigskip

\noindent 2010 {\it Mathematics Subject Classification}:
Primary 05B15; Secondary 68W01. 

\noindent \emph{Keywords:} Collatz conjecture, power series.

\bigskip
\hrule
\bigskip

\noindent (Concerned with sequences.
\seqnum{A002860},
\seqnum{A274171}, and
\seqnum{A274806}.)

\bigskip
\hrule
\bigskip

%\noindent
%Return to
%\htmladdnormallink{Journal of Integer Sequences home page}{https://cs.uwaterloo.ca/journals/JIS/}.
%\vskip .1in

\end{document}